\documentclass[letterpaper, 10 pt, conference]{ieeeconf}  % Comment this line out
\IEEEoverridecommandlockouts                              % This command is only
\usepackage{enumerate}
\usepackage{amsmath}
\usepackage{amsfonts}
\usepackage{graphicx}
\usepackage{algorithm}
\usepackage{algorithmic}
\usepackage{url}
\usepackage{xcolor}
\usepackage{float}

\newtheorem{lemma}{Lemma}

\newtheorem{theorem}{Theorem}

\newtheorem{definition}{Definition}

\newtheorem{remark}{Remark}
\newtheorem{problem}{Problem}
\newcommand {\mat}  [1] {\left[\begin{array}{#1}}
\newcommand {\rix}      {\end{array}\right]}

\title{\LARGE \bf Finding the Nearest Negative Imaginary System with Application to Near-Optimal Controller Design}
\date{}

\author{Mohamed A. Mabrok, IEEE Member% <-this % stops a space
\thanks{Mohamed A. Mabrok with the  Department of Mathematics and Physics, School of Engineering,  Australian University in Kuwait
        {\tt\small m.a.mabrok@gmail.com}}%
}

\begin{document}

\maketitle
\thispagestyle{empty}
\pagestyle{empty}

%%%%%%%%%%%%%%%%%%%%%%%%%%%%%%%%%%%%%%%%%%%%%%%%%%%%%%%%%%%%%%%%%%%%%%%%%%%%%%%%
\begin{abstract}

The negative imaginary (NI) systems theory has attracted interests due to the robustness properties of feedback interconnected NI systems. However, a full output optimal controller-synthesis methodology, for such class of systems, is yet to exist. In order to develop a solution towards  this problem, we first develop a methodology to find the nearest NI system to a non NI system. This later  problem stated as follows: for any linear time invariant (LTI) system defined by the state space matrices  $(A, B, C, D)$, find the nearest NI system, with the state space matrices $(A+\Delta_A,B+\Delta_B,C+\Delta_C,D+\Delta_D)$, such that  the norm of $(\Delta_A,\Delta_B,\Delta_C,\Delta_D)$ is minimized. Then, this methodology will be used to find the nearest optimal controller for a given NI plant.  In other words, for a given NI system, an  optimal control methodology, such as LQG, is used to design an optimal controller that satisfy a particular performance measure. Then, the developed methodology of finding the nearest NI system is used, as a near-optimal control synthesis methodology, to find the nearest NI system to the designed optimal controller. Hence, the synthesized controller satisfy the NI property and therefore guarantee a robust feedback loop with the  negative imaginary system under control.

\end{abstract}

%%%%%%%%%%%%%%%%%%%%%%%%%%%%%%%%%%%%%%%%%%%%%%%%%%%%%%%%%%%%%%%%%%%%%%%%%%%%%%%%
\section{INTRODUCTION}

Dynamical systems with flexibly structured dynamics exists  in many important  engineering systems. For instance, in aerospace systems  \cite{Harigae20,tran2017},  robot manipulators \cite{Wilson2002},  atomic force microscopes (AFMs) \cite{Bhikkaji2009, Mahmood2011}, and other nano-positioning systems \cite{Devasia2007, Diaz2012}). The common  physical structure between these systems is the  involvement of a combination of force actuators and sensors.  These systems are susceptible to changes in environmental variables. For instance, the existence of the highly resonant modes in these flexible systems can affect robustness and stability characteristics \cite{Preumont2011,fanson1990,petersen2010}. Relatively small stimuli or perturbations in the environment, such as changing temperature, can lead to significant impacts in the systems' resonant frequencies. These changes in resonant frequencies result in  changes in the system's phase and gain at a given frequency, which  can easily lead to unacceptable behavior, such as   instability or poor performance in the  system.  Moreover, these dynamical systems, by their nature,  are infinite-dimensional systems \cite{Preumont2011, Morris1993, Ray1978281}, which are hard to model and control. Instead, finite-dimensional models are used as approximation models to be used in controller synthesis \cite{ Preumont2011}.

Negative imaginary (NI) systems theory represent an important class of these  flexibly structure dynamics systems  \cite{lanzon2008,petersen2010,mabrok2012a,Mabrok2015TAC,ferrante2014foundations,Mabrok2015,Mabrok2015d}, In particular, flexible structure dynamics when a force actuator is combined with a collocated, position or acceleration, sensor are known to satisfy the NI property. 

The NI propriety can be  defined, in the case of the single-input single-output (SISO),  based on  the properties of the imaginary part of the complex function of the   frequency response $G(j\omega)$. 
 This is,   
 \begin{equation*}
     j\left( G(j\omega )-G(j\omega )^{\ast}\right) \geq 0,
 \end{equation*}
to be satisfied for all $\omega\in(0,\infty)$. A formal definition of NI system is given in the next section . 

\begin{figure}
  \centering\includegraphics[width=7cm]{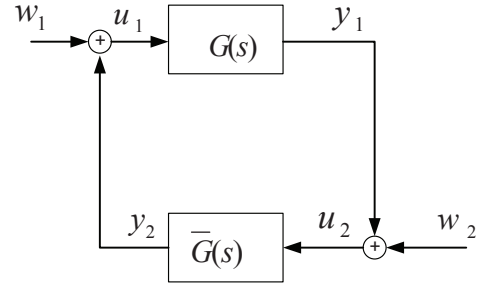}\\
  \caption{A negative-imaginary feedback control system. 
}\label{conn:NI:SNI}
\end{figure}

A fundamental result of NI systems theory is concerning the feedback of interconnected NI systems. That is, the closed feedback interconnection between an NI system and strictly negative imaginary, as shown in Fig. \ref{conn:NI:SNI}, is stable under a given DC gain condition. This result made the NI theory attractive to be used in controller design, i.e., if the system  under control satisfy the NI property, one can design an NI controller so that the robust stability is guaranteed for free. This led to several attempts to synthesis an NI controller that preserve certain performance measures \cite{lanzon2009,petersen2010,xiong2016,Mabrok2015_IET,mabrok2019}.   However, a full output optimal controller for an NI system is yet to exist. One aim of this  paper is  to present  a near-optimal control synthesis methodology for negative imaginary systems. The synthesized  controller satisfies the NI property, and therefore, guarantee a robust feedback loop with the  negative imaginary system under control. 

In order to develop the near-optimal control  methodology, we first develop a method to find the nearest NI system for any non-NI system. In other words, given an LTI system, 
minimize the Frobenius norm of $(\Delta_A,\Delta_B,\Delta_C,\Delta_D)$ such that $(A+\Delta_A,B+\Delta_B,C+\Delta_C,D+\Delta_D)$ satisfies the NI property. The problem of finding the nearest negative imaginary system is motivated by a similar problem of finding the nearest positive real system (passive system) presented in \cite{gillis2018finding,fazzi2021finding,SchT07,WanZKPW10,VoiB11}.  In  \cite{SchT07}, several  assumptions, such as  $D+D^T$ to be  non-singular are made. Also, they  restrict the  perturbation on the matrix $C$ only. The methods in~\cite{WanZKPW10} and~\cite{VoiB11} they perturb, both,  $B$ and/or $C$. In~\cite{BruS13}, perturbations in all system matrices. Similar results for NI system were developed in    \cite{mabrok2013enforcing} where an algorithm was developed for enforcing negative imaginary property in case of any violation during system identifications. This assumes that the underlying dynamics ought to belong to negative imaginary  system class. The  method is based on  the spectral properties of Hamiltonian matrices. 

In this paper, we use the result developed in \cite{gillis2018finding} to develop similar results of finding the  nearest negative imaginary system. One of the main advantages of this method over the other perturbation methods is that there are no assumptions in the system. Also, in the positive real case, it allows for perturbations of all system matrices $(A,B,C,D)$.

The developed method of finding the nearest NI system will be used in order to find a near-optimal NI controller for a given NI plant. We first use a regular optimal control methodology such as LQG to design an optimal controller for a given NI system. Then, we employ the developed method of the  nearest NI dynamical system to find the NI optimal controller.

As discussed, there are two main contributions in this paper, which can be summarized as follows:  
\begin{enumerate}
    \item  First, the paper introduces a methodology for finding the nearest negative imaginary system for a non-negative imaginary system. The proposed method is based on the  Port-Hamiltonian formulation of the negative imaginary systems. 
    \item Second, these result of finding the nearest negative imaginary system  will be used in order to find the  nearest NI optimal controller for a given negative imaginary plant. The synthesized controller satisfies a near optimal performance with the negative imaginary robustness property.  
\end{enumerate}

% \textcolor{red}{TODO}
% This paper is further organized as follows: Section \ref{sec:pre} presents the  preliminaries that are used thought the paper.  Subsection \ref{phsys} presents the Port-Hamiltonian formulation of the negative imaginary systems. The problem  formulation of finding the nearest negative imaginary system is dissuaded in \ref{sec:theproblem}. An algorithmic for nearest negative imaginary system problem is presented in \ref{algosol}. The optimal control design problem is discussed in \ref{sec:LQG}.

\section{Preliminaries}
In this section, we present the definitions and lemmas of the negative imaginary systems. Also, this section presents the Port-Hamiltonian formulation of the negative imaginary systems.
 \label{sec:pre}
\subsection{Negative imaginary systems }

In this section, we recall  the definition of linear time invariant (LTI) negative imaginary systems  as given in \cite{Mabrok2015TAC}.

Consider the following LTI system defined in \eqref{eq:ltisys}
\begin{align}
\label{eq:ltisys}
\begin{split}
\dot{x}(t)= & Ax(t)+Bu(t), \\
y(t) = & \; C x(t) + Du(t),
\end{split}
\end{align}
where, $A$ represents the dynamical matrix with dominion belongs to  $\mathbb{R}^{n \times n}$, $B$ is the input matrix $ \in \mathbb{R}^{n \times m}$, $C$ represents the output matrix 
$\in \mathbb{R}^{m \times n},$ and  $D \in \mathbb{R}^{m \times m}.$ This yields a transfer function matrix $G(s)=C(sI-A)^{-1} B+D$.
The transfer function matrix $G(s)$ is said to be strictly proper if
$G(\infty)=D=0$. 
The notation $
\begin{bmatrix}
\begin{array}{c|c}
A & B \\ \hline C & D
\end{array}
\end{bmatrix}$ will be used  to denote the state space realization
\eqref{eq:ltisys}.

The  NI system is formally defined in \cite{Mabrok2015TAC}.

\begin{definition}\cite{Mabrok2015TAC}\label{Def:NI}
A square transfer function matrix $G(s)$ is NI if the following conditions are satisfied:
\begin{enumerate}
\item $G(s)$ has no pole in $Re[s]>0$.
\item For all $\omega >0$ such that $s=j\omega$ is not a pole of $G(s)$,
\begin{equation}\label{eq:NI:def}
    j\left( G(j\omega )-G(j\omega )^{\ast }\right) \geq 0.
\end{equation}
\item If $s=j\omega _{0}$ with $\omega _{0}>0$ is a pole of $G(s)$, then it is a simple pole and the residue matrix $K=\underset{%
s\longrightarrow j\omega _{0}}{\lim }(s-j\omega _{0})jG(s)$ is
Hermitian and  positive semidefinite.
 \item If $s=0$ is a pole of $G(s)$, then
$\underset{s\longrightarrow 0}{\lim }s^{k}G(s)=0$ for all $k\geq3$
and $\underset{s\longrightarrow 0}{\lim }s^{2}G(s)$ is Hermitian and
positive semidefinite.
\end{enumerate}
\end{definition}
%\end{definition}

\begin{definition}\cite{xiong21010jor}
A square transfer function matrix $G(s)$ is SNI if  the
following conditions are satisfied:
\begin{enumerate}
\item ${G}(s)$ has no pole in $Re[s]\geq0$.
\item For all $\omega >0$, $j\left( {G}(j\omega )-{G}(j\omega )^{\ast }\right) > 0$.
\end{enumerate}
\end{definition}

Here, we present an NI lemma.
\begin{lemma}\label{NI-lemma}
 Let $
\begin{bmatrix}
\begin{array}{c|c}
A & B \\ \hline C & D
\end{array}
\end{bmatrix}$ defining the system (\ref{eq:xdotn})-(\ref{eq:yn}) be a minimal realization of the transfer function matrix  $%
G(s)$. Then, $G(s)$ is NI if and only if $D=D^T$ and there exist
matrices $P=P^{T}\geq 0$,
 $W\in \mathbb{R}^{m \times m}$, and $L\in \mathbb{R}^{m \times n}$
 such that the following linear matrix inequality (LMI) is
satisfied:
\begin{small}
\begin{align}\label{LMI:PR}
\begin{bmatrix}
PA+A^{T}P & PB-A^{T}C^{T} \\
B^{T}P-CA & -(CB+B^{T}C^{T})%
\end{bmatrix}%
 =
\begin{bmatrix}
-L^{T}L & -L^{T}W \\
-W^{T}L & -W^{T}W%
\end{bmatrix}%
\leq 0.
\end{align}
\end{small}
\end{lemma}

A state-space characterization of NI systems in terms of a pair of
linear matrix inequalities (LMIs) has been given in
\cite{lanzon2008}. This result also
generalized in \cite{xiong21010jor} to include poles on the
imaginary axis except at the origin.

\begin{lemma}\label{ch:li:NI-lemma}(See \cite{xiong21010jor})
Let $
\begin{bmatrix}
\begin{array}{c|c}
A & B \\ \hline C & D
\end{array}
\end{bmatrix}$
be a minimal state space realization of a transfer function matrix
$G(s)$. Then $G(s)$ is NI if and only if $det(A)\neq0$, $D=D^{T}$
and there exists a real matrix $Y>0$ such that
\begin{equation}
AY+YA^{*}\leq0 \ \text{and}\ B=-AYC^{*}.
\end{equation}
\end{lemma}

\begin{lemma}\label{ch5:NI-lemma}
 Let $
\begin{bmatrix}
\begin{array}{c|c}
A & B \\ \hline C & D
\end{array}
\end{bmatrix}$ be a minimal realization of the transfer function matrix  $%
G(s)$ for the system in \eqref{eq:ltisys}. Then,
$G(s)$ is NI if and only if $D=D^T$ and there exists a matrix
$P=P^{T}\geq 0$
 %$W\in \mathbb{R}^{m \times m}$, and $L\in \mathbb{R}^{m \times n}$
 such that the following LMI is
satisfied:
\begin{align}\label{ch5:LMI:PR}
\begin{bmatrix}
PA+A^{T}P & PB-A^{T}C^{T} \\
B^{T}P-CA & -(CB+B^{T}C^{T})%
\end{bmatrix}%
%&=%
%\begin{bmatrix}
%-L^{T}L & -L^{T}W \\ \nonumber
%-W^{T}L & -W^{T}W%
%\end{bmatrix}\\
\leq 0.
\end{align}
Furthermore, if $G(s)$ is SNI, then $\det(A)\neq0$ and there exists
a matrix $P>0$ such that \eqref{ch5:LMI:PR} holds.
\end{lemma}

One of the important results in the NI system's theory is the robustness property that emerge in the case of  positive feedback interconnection between an NI and an SNI system as shown in Fig.
\ref{conn:NI:SNI}. The following theorem form \cite{petersen2010,lanzon2008} states this results: 
\begin{theorem}
\label{ch:li:th:ni-con} Consider an NI transfer function matrix
$G(s)$ with no poles at the origin and an SNI transfer function
matrix $\bar{G}(s)$, and suppose that $G(\infty)\bar{G}(\infty)= 0$
and $\bar{G}(\infty)\geq  0$. Then, the positive-feedback
interconnection  (see Fig. \ref{conn:NI:SNI}) of $G(s)$ and
$\bar{G}(s)$  is internally stable if and only if
$\lambda_{max}(G(0)\bar{G}(0))<1$.
\end{theorem}

The above Theorem \ref{ch:li:th:ni-con} characterizes  the conditions of the  
stability of the feedback interconnection of two  NI and SNI systems
through the phase stabilization.

\subsection{Port-Hamiltonian formulation of negative imaginary systems} \label{phsys}
To establish the results in this paper, the formulation of the negative imaginary system is used. Port-Hamiltonian formulation of negative imaginary systems was developed in \cite{van2011positive,van2016interconnections}. 
The following lemma characterize the NI systems as a Port-Hamiltonian system. 

\begin{lemma}\label{lem:NI:PH}
The system given in \eqref{eq:ltisys} has negative imaginary transfer function if and only if it can be written as

\begin{align}
\label{eq:phsystem}
\begin{split}
\dot{x}(t) = & \; (J-R)Q(x(t) - C^{T}u(t)), \\
y(t) = & \; C x(t) + Du(t),
\end{split}
\end{align}

for some matrices $Q, J, R$, where, 
\begin{equation}
    Q=Q^T>0,\ \ \ J=-J^T, \  \  \ R=R^T\geq 0.
\end{equation}
\end{lemma}

The proof is the same as proposition (2.4) in \cite{van2011positive}.

The above lemma  be our main tool to reformulate the nearest negative imaginary  system problem. The next sections will focus on defining the problem formulation.

%==================
%==================
%==================

\section{Main results}
 In this section, we present the main results in this paper. 

\subsection{Nearest negative imaginary system problems} \label{sec:theproblem}

As indicated in the introduction,  the problem of finding the nearest negative imaginary system is similar  to the problem of finding the nearest positive real system (passive system) presented in \cite{gillis2018finding}, where the  port-Hamiltonian formulation  is used.

Let us  now define the nearest negative imaginary  system problems under consideration:

\begin{problem}\label{prob_g}
Suppose an LTI system with the following  state space representation $
\begin{bmatrix}
\begin{array}{c|c}
A & B \\ \hline C & D
\end{array}
\end{bmatrix}$, find the nearest (the closest)  system $
\begin{bmatrix}
\begin{array}{c|c}
 \tilde A & \tilde B \\ \hline \tilde C & \tilde D
\end{array}
\end{bmatrix},$
%$(\tilde A,\tilde B,\tilde C, \tilde D)$ to $(A,B,C,D)$,
such that, 
\begin{equation*}\label{def_F}
\inf_{(\tilde A,\tilde B,\tilde C, \tilde D) }
\mathcal{F}(\tilde A,\tilde B,\tilde C,\tilde D),
\end{equation*}
where,
\begin{align}\label{eq:def_F}
\mathcal{F}(\tilde A,\tilde B,\tilde C,\tilde D) =&
{\|A-\tilde A\|}_F^2+{\|B-\tilde B\|}_F^2 \nonumber\\
&+{\|C-\tilde C\|}_F^2+{\|D-\tilde D\|}_F^2.
\end{align}
\end{problem}

The following definition is based on  Lemma \ref{lem:NI:PH} and compares the system descried in  \eqref{eq:phsystem} with the LTI system given in \eqref{eq:ltisys}. 

\begin{definition}{\rm
A system $(A,B,C,D)$ is said to admit a port-Hamiltonian form if
there exists a  system as defined in~\eqref{eq:phsystem} such that
\[
A=(J-R)Q, \text{ and}, B = -(J-R)C^T
\]
}
\end{definition}

Based on the above definition, the problem  given in \eqref{prob_g} can be reduced to the following problem: 

\begin{problem}\label{prob_g2}
% For a given system $(A,B,C,D)$, find the nearest system $(\tilde A,\tilde B,\tilde C, \tilde D)$ to $(A,B,C,D)$, such that, 
Suppose an LTI system with the following  state space representation $
\begin{bmatrix}
\begin{array}{c|c}
A & B \\ \hline C & D
\end{array}
\end{bmatrix}$, find the nearest (the closest)  system $
\begin{bmatrix}
\begin{array}{c|c}
 \tilde A & \tilde B \\ \hline \tilde C & \tilde D
\end{array}
\end{bmatrix},$
%$(\tilde A,\tilde B,\tilde C, \tilde D)$ to $(A,B,C,D)$,
such that, 
\begin{equation*}\label{def_F}
\inf_{(\tilde A,\tilde B,\tilde C, \tilde D) }
\mathcal{F}(\tilde A,\tilde B,\tilde C,\tilde D),
\end{equation*}
where,
\begin{equation}\label{eq:def_F}
\mathcal{F}(\tilde A,\tilde B) =
{\|A-(J-R)Q\|}_F^2+{\|B-(R-J)C^T\|}_F^2,
\end{equation}
where, $ \tilde A=(J-R)Q, \ \  \tilde B=-(J-R)C^T$. 
\end{problem}

\subsection{Algorithmic for nearest negative imaginary system problem} \label{algosol}

This section proposes an algorithm to solve the problem dissuaded in the above section.

The problem \eqref{prob_g2} can be written  as follows
\begin{align}
 \inf_{J,R,Q} & %\mathcal{F}((J-R)Q,F-P,(F+P)^TQ,S+N,I_n)
 {\|A-(J-R)Q\|}_F^2 + {\|B-(R-J)Q\|}_F^2,
\nonumber \\
& \text{ such that } J^T=-J, Q=Q^T >0, \text{and } R^T=R\geq 0. \label{def:dist_sph_1}
\end{align}

The projected gradient method (FGM) presented in \cite{GilS16} and in \cite{gillis2018finding} is used to solve the problem in ~\eqref{def:dist_sph_1}.

As indicated in \cite{gillis2018finding}, the  projected gradient method is much faster  and hence better to use  compared to the standard projected gradient method~\cite{nes83}. 

The steps can be summarized as follows: 
\begin{itemize}
\item Compute the gradient as follows: 
\begin{small}
\begin{align*}
    &{\|A-(J-R)Q\|}_F^2 \Rightarrow \nabla_{J_1 }= -2Q{\|A-(J-R)Q\|},\\
    &{\|A-(J-R)Q\|}_F^2 \Rightarrow \nabla_{R_1} = 2Q{\|A-(J-R)Q\|},\\
   & {\|A-(R-J)C^T\|}_F^2 \Rightarrow \nabla_{J_2} = 2C^T{\|A-(R-J)C^T\|},\\
    &{\|A-(R-J)C^T\|}_F^2 \Rightarrow 
    \nabla_{R_2} = -2C^T{\|A-(R-J)C^T\|},\\
        & {\|A-(J-R)Q\|}_F^2 \Rightarrow 
         \nabla_Q =- 2(J-R){\|A-(J-R)Q\|},
\end{align*}
\end{small}
 or simply, for a given term in the objective function, $f(X) = {\|AX-B\|}_F^2$
the  gradient is $\nabla_Y f(Y) = 2A^T(AY-B)$.
\item Project onto the feasible set of matrices $Q, R$ that satisfy both conditions,   $ Q=Q^T >0, \text{and } R^T=R\geq 0$.
\end{itemize}
  The FGM Algorithm, which presented in  \cite{gillis2018finding}, is used to compute the matrices $Q, R$.

Similar to the  implementation in \cite{gillis2018finding}, positive weights $w_i$ were added to the objective function terms in order to give opportunity for  a different importance of each term if needed. Therefore, the objective function can be written as follows: 
\begin{small}
\begin{equation*}\label{eq:weighted_obj}
\mathcal{F}(\tilde A,\tilde B) = w_1 {\|A-(J-R)Q\|}_F^2 + w_2 {\|B-(R-J)C^T\|}_F^2.
\end{equation*}
\end{small}
Parameter settings in  our implementation are  similar to the  parameter settings that was used in \cite{gillis2018finding}. For instance, the  step length is calculated as  $\gamma = 1/L$ where $L={\|Q\|}_2^2$. Moreover, in the initialization step, two different initializations were used.  

%\begin{itemize}
\textbf{The first initialization},
\[ Q=I_n, \ \ \ \ \
J = \big(A-A^T\big)/2, \ \ \ \
R =  P_{>0 }\big((-A-A^T)/2\big),
\]
where the  notation $\mathcal P_{>0 }(X)$ stands for the projection of a matrix $X$
on the cone of positive semi-definite matrices.

 \textbf{LMI-based initialization}: Since the given  system is not an NI system,  the LMI given in \eqref{ch5:LMI:PR} has no solution. However, a solution $P$ to nearby LMIs should be a good initialization to the matrix $Q$.
We propose the following to relax the LMIs~\eqref{ch5:LMI:PR}:
\begin{align}
&\min_{\delta, P}  \quad \delta^2 \nonumber \\
& \ \ \text{ such that }  \quad \nonumber  \\
&\begin{bmatrix}
-PA-A^{T}P & -PB+A^{T}C^{T} \\
-B^{T}P+CA & (CB+B^{T}C^{T})%
\end{bmatrix} + \delta I_{n+m} \geq 0,  \label{initLMIs} 
\end{align}

\subsection{Optimal control design}\label{sec:LQG}
In this section, the nearest NI problem, which was presented in the previous  subsection,  will be used in order to design a  near-optimal controller for a given NI plant. 

Suppose that we want to design a controller for a given  NI plant, $G(s)$, with the state space representation given in \eqref{eq:ltisys}. Suppose also that we decided to use any standard control synthesis  methodology such as LQG or $H_\infty$ to design a controller that satisfy a particular performance  measure. It is unlikely that the designed controller will satisfy the NI property and therefore, a robustness property will not be guaranteed. Hence, we can use the  nearest NI problem, which was presented in the previous  subsection, to find the nearest NI controller to the designed one. 
The following steps summarize the NI-control design, assuming that an LGQ is used in the design. 
\begin{itemize}
    \item Given an LTI NI plant  in the form \eqref{eq:ltisys}, with the transfer function matrix $G(s)=C(sI-A)^{-1} B+D$.
    \item Design a linear quadratic Gaussian (LQG) controller $K(s)= C_c(sI-A_c)^{-1} B_c+D_c$, which minimizes the following  cost function:
    \begin{equation}
        J(u)=\int_{0}^{\infty}\left\{x^{T} Q_c x+2 x^{T} N_c u+u^{T} R_c u\right\} d t.
    \end{equation}
    \item Use the methodology presented in this paper to find the nearest NI controller $\bar{G}(s)=\bar{C}(sI-\bar{A})^{-1} \bar{B}+\bar{D},$ to the designed LQG controller $K(s)$. 
\end{itemize}
The new modified controller $\bar{G}(s)$ is a near-optimal controller that satisfy the NI property.

\begin{remark}
The DC gain condition  $\lambda_{max}(G(0)\bar{G}(0))<1, $ can be included in the optimization process of finding the nearest NI controller. The DC gain of the NI controller can be calculated as follows: 
\begin{align*}
   \bar{G}(0) =&-\bar{C}(\bar{A})^{-1} \bar{B}+\bar{D} ,\\
              =& \ \bar{C}((J-R)Q)^{-1} (J-R)\bar{C}^T+\bar{D}, \\
              =& \ \bar{C}Q^{-1}\bar{C}^T+\bar{D}.
\end{align*}
In the iterations of finding the matrix $Q$, particularly, in the projection iteration, the matrix  $Q$ is scaled to satisfy the DC gain condition. The scaling factor that preserve the DC gain condition is:
\begin{align*}
    Q_{new} = \alpha Q_{old}
\end{align*}
where in the single-input single-output case,
\begin{equation*}
    \alpha = (CQ^{-1}_{old}C^T+D)G(0)+\epsilon,
\end{equation*}
 with a small $\epsilon >0$.
\end{remark}

\section{Example} \label{numexp}
In this section, we present an example to illustrate the design approach presented in this paper. 

It is well known that mechanical structures with colocated force actuators and position sensors yield negative imaginary systems \cite{petersen2010}. Naturally, these systems are infinite dimension systems, whereas their models are not. This make the control design for such systems challenging. Particularly, in the case where the synthesis methodology  do not take into account the robustness issue. Therefore, our method shows a big advantage over optimal control methodologies. 

\begin{figure}[h!]\label{PlantBodeN2}
\centering
\includegraphics[width=0.5\textwidth]{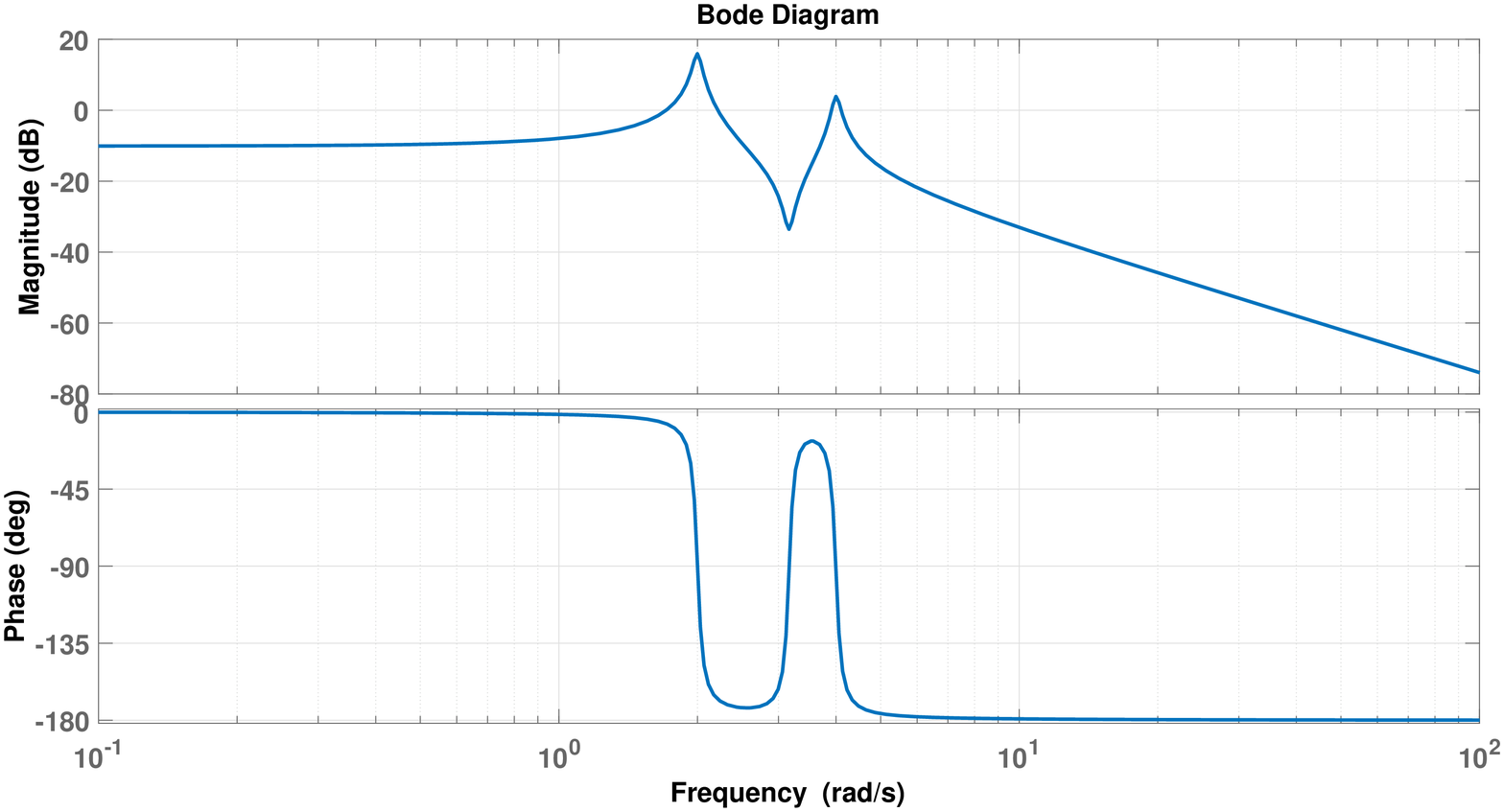}
\caption{Bode plot of NI system given in \eqref{lti:SOSG}, where N=2}
\end{figure}

To illustrate this fact, consider the following lightly damped flexible structure LTI second-order system with a colocated force actuation
and position measurement with the following structure: 
\begin{align}\label{lti:SOSG}
    G(s) = \sum_{n=1}^{N}\frac{1}{s^{2}+2 \zeta_n \omega_{n} s+\omega_{n}^{2}},
   %\frac{1}{ s^2 + 0.16 s + 16}
\end{align}
where $\omega_{n} $ is  the natural frequency and $\zeta_n$ the damping factor. Suppose that we want to design an LQG controller for the system given in \eqref{lti:SOSG}. Since this model represents an infinite dimension system, a finite model is chosen to design the controller. We chose $N=2$ with $\omega_{1}= 2, \omega_{2}= 4$ and $\zeta_{1}=\zeta_{2}= 0.02$ for the model parameters. This implies that the model gives the transfer function given in Fig. \ref{PlantBodeN2}. 

\begin{figure}[H]\label{LQGBode}
\centering
\includegraphics[width=0.5\textwidth]{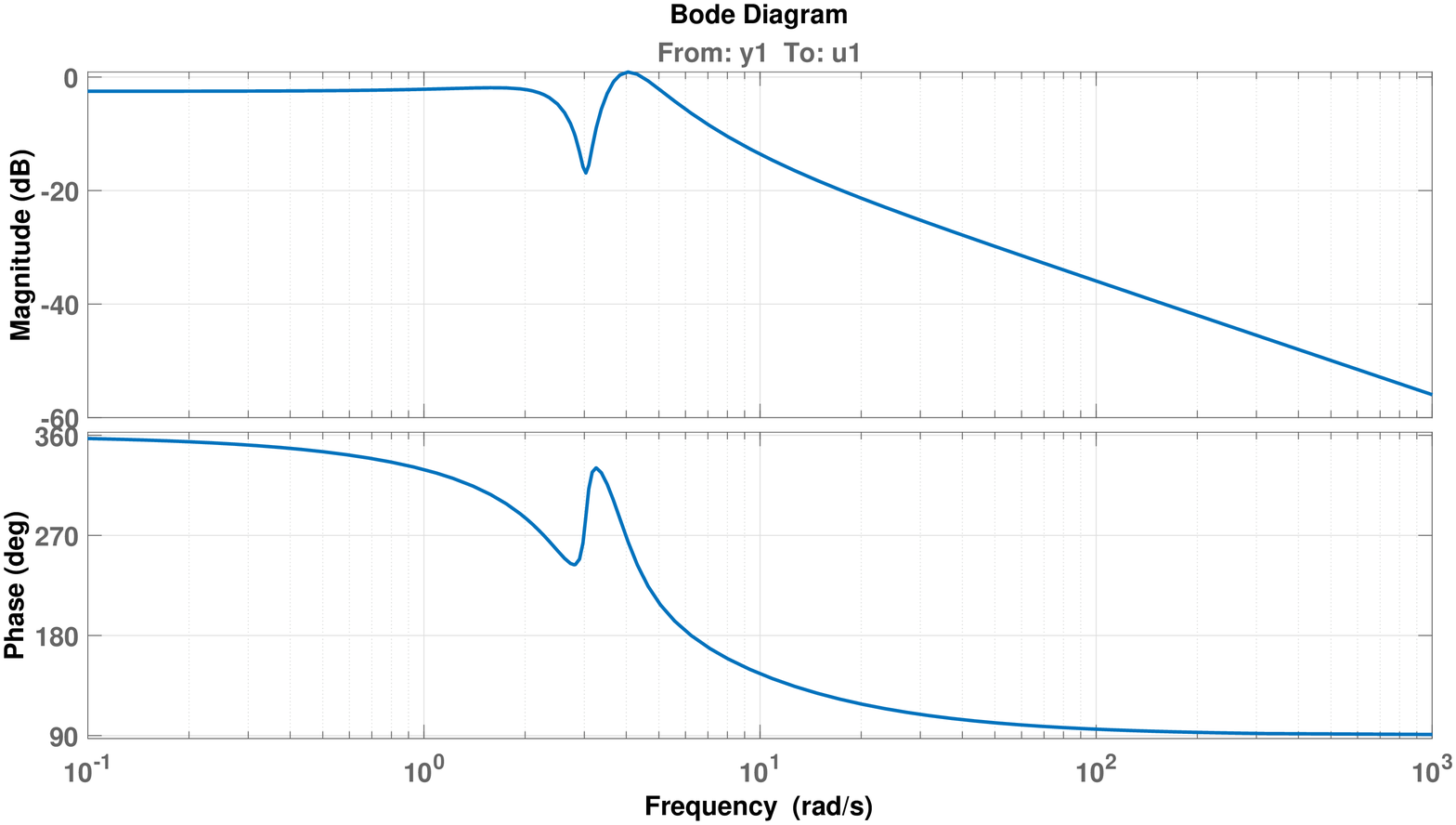}
\caption{Bode plot of designed LQG controller  given in \eqref{LQG1}}
\end{figure}

With an appropriate LQG parameters, the  controller is given  as follows:
\begin{align}\label{LQG1}
   LQG(s)=\frac{ -1.593 s^3 + 9.84 s^2 - 12.58 s + 93.76}{s^4 + 3.847 s^3 + 26.66 s^2 + 46.86 s + 125.1}. 
\end{align}
The bode plot of the designed LGQ controller as given in \ref{LQGBode} shows that it is not an NI controller, since the phase is not in the $(0,-\pi)$. 

\begin{figure}[h!]\label{NILQGBode}
\centering
\includegraphics[width=0.5\textwidth]{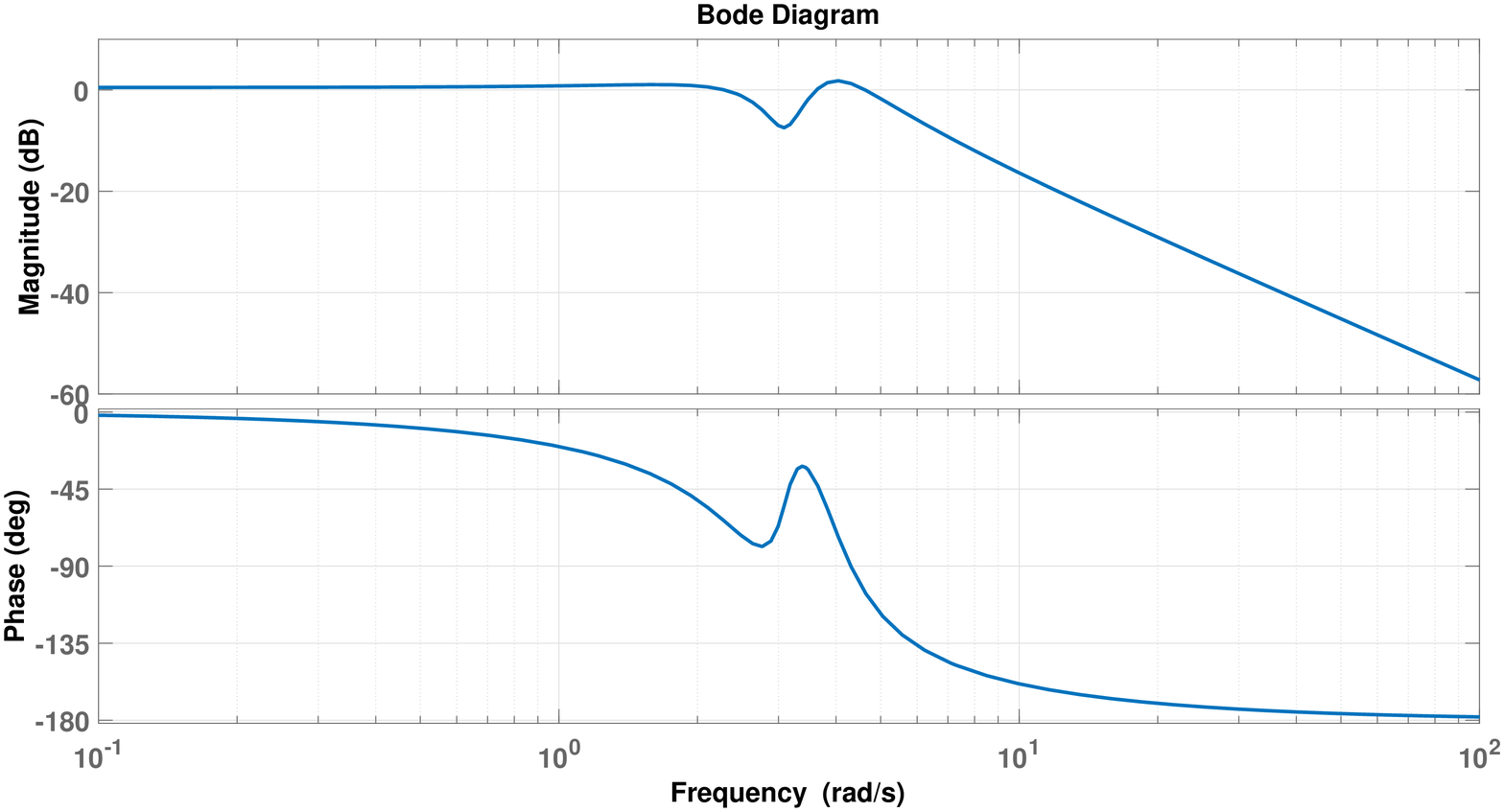}
\caption{Bode plot of designed LQG controller  given in \eqref{LQG1}}
\end{figure}

Now, applying our method of finding the nearest NI controller to the LQG controller given in  \eqref{LQG1}, we get the following transfer function,

\begin{align}\label{NILQG1}
  NILQG(s)= \frac{ 13.75 s^2 + 6.77 s + 132.5}
  {s^4 + 3.847 s^3 + 26.66 s^2 + 46.86 s + 125.1}
\end{align}
The bode plot in Fig. \ref{NILQGBode} of the  controller given in \eqref{NILQG1} shows that it satisfy the NI property.

Applying FGM on with the standard initialization,
This gives a nearby standard NI system with error
\begin{small}
\[
{\|A-\hat A\|}_F^2 +{\|B-\hat B\|}_F^2+ {\|C-\hat C\|}_F^2+{\|D-\hat D\|}_F^2 = 0.6430.
\]
\end{small}
In terms of relative error for each matrix, we have
{\small{\begin{align*}
&\frac{{\|A-\hat A\|}_F}{{\|A\|}_F} = 5.5917e^{-18}  \%,~
\frac{{\|B-\hat B\|}_F}{{\|B\|}_F} = 0.0631  \%,~\\
&\frac{{\|C-\hat C\|}_F}{{\|C\|}_F} = 0  \%,~ 
\frac{{\|D-\hat D\|}_F}{{\|D\|}_F} = 0  \%.
\end{align*}}}

\begin{figure}[H]\label{StepBothLQGNI}
\centering
\includegraphics[width=0.5\textwidth]{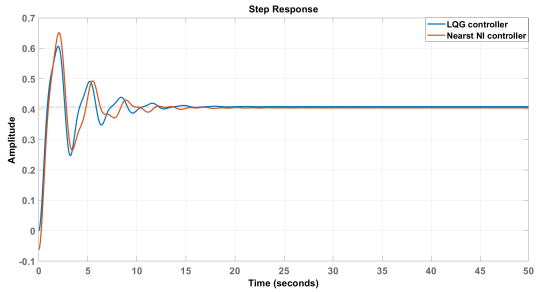}
\caption{Bode plot of NI system given in \eqref{lti:SOSG}, where N=2}
\end{figure}

The step response of the  closed feedback interconnection of the  plant given in \eqref{lti:SOSG} and both the  designed LQG \eqref{LQG1} controller and the nearest NI controller \eqref{NILQG1} are  given in Fig. \ref{StepBothLQGNI}. It is clear that the response is very similar.  

\begin{figure}[H]\label{BodeBothplants}
\centering
\includegraphics[width=0.5\textwidth]{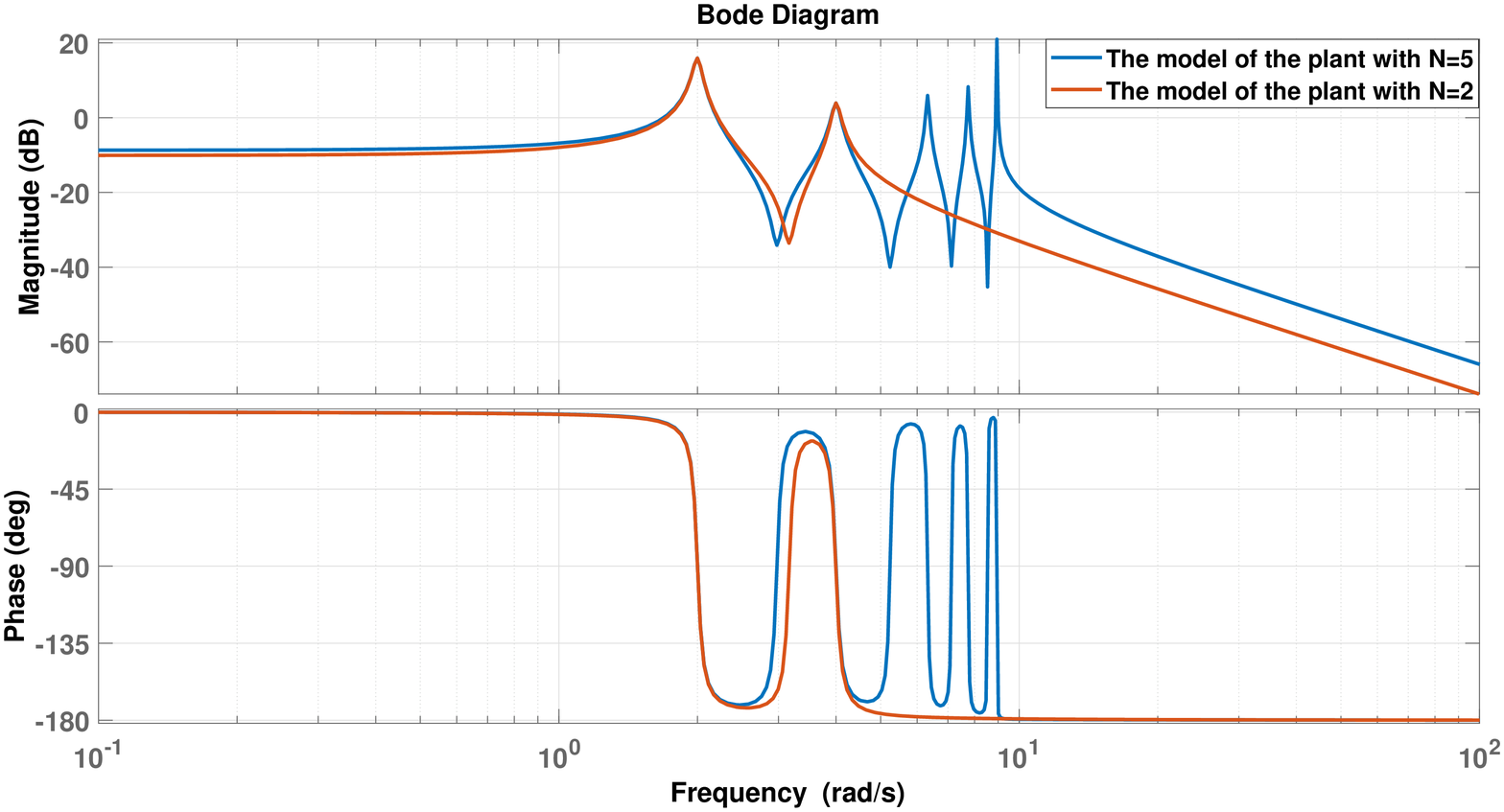}
\caption{Bode plot of the plant model  given in \eqref{lti:SOSG}, with different numbers modes}
\end{figure}

The more interesting part of this example is when we add more non-modeled modes to the plant, i.e., $N= 5$ as shown in Fig. \ref{BodeBothplants} for instance. This means that we include some of the un-modulated dynamics in the plant, which was regarded as uncertainty. For instance, suppose that the plant  

\begin{figure}[H]\label{BodeBothplantsunstable}
\centering
\includegraphics[width=0.5\textwidth]{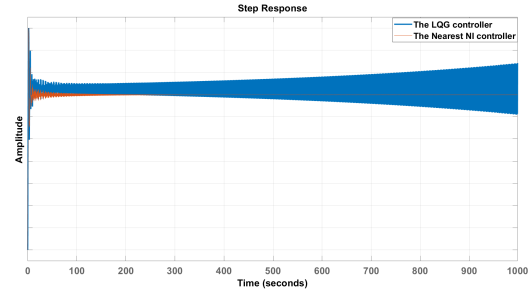}
\caption{The step response of the designed LGQ \eqref{LQG1} and the nearest NI  controller  plant model  given in \eqref{lti:SOSG}, with  N=5}
\end{figure}

As shown in Fig. \ref{BodeBothplantsunstable}, the designed LGQ \eqref{LQG1} will become unstable if we considered the five-mode plant. However, the nearest NI  controller still stabilize the system with an acceptable performance. This is due to the negative imaginary property of both, the controller and the plant.  

% \bibliographystyle{ieeetr}
% \bibliography{ref.bib}

\end{document}